\documentclass[12pt]{amsart}

\usepackage{amssymb, geometry}

\begin{document}
\baselineskip=18pt
\setcounter{page}{1}
    
\newtheorem{conjecture}{Conjecture\!\!}
\newtheorem{remark}{Remark\!\!}
\newtheorem{theorem}{Theorem}
\newtheorem{lemma}{Lemma\!\!}
\newtheorem{proposition}{Proposition\!\!}

\renewcommand{\theconjecture}{}
\renewcommand{\theremark}{}
\renewcommand{\theproposition}{}
\renewcommand{\thelemma}{}

\def\a{\alpha}
\def\aa{{\bf a}}
\def\bb{{\bf b}}
\def\cc{{\bf c}}
\def\dd{{\bf d}}
\def\b{\beta}
\def\B{{\bf B}} 
\def\C{{\bf C}} 
\def\K{{\bf K}}
\def\BB{{\mathcal{B}}} 
\def\DD{{\mathcal{D}}} 
\def\GG{{\mathcal{G}}} 
\def\KK{{\mathcal{K}}} 
\def\LL{{\mathcal{L}}} 
\def\SS{{\mathcal{S}}}
\def\UU{{\mathcal{U}}}
\def\ca{c_{\a}}
\def\ka{\kappa_{\a}}
\def\coa{c_{\a, 0}}
\def\cua{c_{\a, u}}
\def\cL{{\mathcal{L}}} 
\def\cD{{\mathcal{D}}} 
\def\Ea{E_\a}
\def\eps{{\varepsilon}} 
\def\esp{{\mathbb{E}}} 
\def\Ga{{\Gamma}} 
\def\G{{\bf \Gamma}} 
\def\e{{\rm e}}
\def\i{{\rm i}}
\def\L{{\bf L}}
\def\lbd{\lambda}
\def\lacc{\left\{}
\def\lcr{\left[}
\def\lpa{\left(}
\def\lva{\left|}
\def\M{{\bf M}}
\def\Nn{{\bf N}}
\def\prst{{\leq_{st}}}
\def\prost{{\prec_{st}}}
\def\prcvx{{\prec_{cx}}}
\def\Rr{{\bf R}}
\def\CC{{\mathbb{C}}}
\def\NN{{\mathbb{N}}} 
\def\QQ{{\mathbb{Q}}} 
\def\pb{{\mathbb{P}}}
\def\racc{\right\}}
\def\rl{{\mathbb{R}}}
\def\rpa{\right)}
\def\rcr{\right]}
\def\rva{\right|}
\def\Tt{{\bf \Theta}}
\def\Ttt{{\tilde \Tt}}
\def\W{{\bf W}}
\def\X{{\bf X}}
\def\XX{{\mathcal X}}
\def\YY{{\mathcal Y}}
\def\U{{\bf U}}
\def\V{{\bf V}}
\def\Un{{\bf 1}}
\def\Y{{\bf Y}}
\def\Z{{\bf Z}}
\def\ZZ{{\mathbb{Z}}}
\def\A{{\bf A}}
\def\AA{{\mathcal A}}
\def\hAA{{\hat \AA}}
\def\hL{{\hat L}}
\def\hT{{\hat T}}

\def\claw{\stackrel{d}{\longrightarrow}}
\def\elaw{\stackrel{d}{=}}
\def\pslaw{\stackrel{a.s.}{\longrightarrow}}
\def\qed{\hfill$\square$}

\newcommand*\pFqskip{8mu}
\catcode`,\active
\newcommand*\pFq{\begingroup
        \catcode`\,\active
        \def ,{\mskip\pFqskip\relax}%
        \dopFq
}
\catcode`\,12
\def\dopFq#1#2#3#4#5{%
        {}_{#1}F_{#2}\biggl[\genfrac..{0pt}{}{#3}{#4};#5\biggr]%
        \endgroup
}


\newcommand{\R}{\mathbb{R}}
\newcommand{\N}{\mathbb{N}}
\newcommand{\bP}{\mathbb{P}}    
\newcommand{\bE}{\mathbb{E}}    
\def\ii{{\rm i}}
\def\S{\mathbf{S}}
\def\F{\mathbf{T}}
\def\W{\mathbf{W}}

\title{On some new moments of Gamma type}

\author[T.~Kadankova]{Tetyana Kadankova}

\address{Vakgroep Wiskunde, Vrije Universiteit Brussel, Gebauw G, Pleinlaan 2, 1050 Elsene, Belgium. {\em Email}: {\tt tetyana.kadankova@vub.be}}

\author[T.~Simon]{Thomas Simon}

\address{Laboratoire Paul Painlev\'e, Universit\'e de Lille, Cit\'e Scientifique, 59655 Villeneuve d'Ascq Cedex, France. {\em Email}: {\tt thomas.simon@univ-lille.fr}}

\author[M.~Wang]{Min Wang}

\address{Jingzhai 301, Yau Mathematical Sciences Center, Tsinghua University, Haidian District, Beijing 100084, China. {\em Email}: {\tt minwangmath@tsinghua.edu.cn}}



\begin{abstract} We investigate certain positive random variables having moments of Gamma type. Some necessary and some sufficient conditions are given for their existence. In particular, we observe that the Weber-Schafheitlin formula for the Bessel function makes it possible to construct non-trivial moments of Gamma type having a signed spectral measure.
 
\end{abstract}

\maketitle

\section{Introduction and statement of the results}

\subsection{The general setting}\label{general} Introduce the notation
$$({\bf a})_s\; =\; (a_1)_s\,\times\,\cdots\,\times\,(a_n)_s\; =\;\prod_{i=1}^n \frac{\Ga(a_i +s)}{\Ga(a_i)}$$
for ${\bf a} = \{a_1, \ldots, a_n\}$ with $0 < a_1\le \ldots\le a_n$ and $s > -a_1,$ where we have used the standard definition of the Pochhammer symbol $(x)_s$ built on the Gamma function $\Gamma(z), z >0.$ When $n = 0$ that is ${\bf a}$ is empty, we set $\min({\bf a}) = \infty$ and $({\bf a})_s = 1.$ In \cite{CL}, the authors introduce the notation 
$${\rm D}\lcr\begin{array}{cc} {\bf a} & {\bf b}\\ {\bf c} & {\bf d}\end{array}\rcr$$ 
for the distribution of the positive random variable $X,$ if it exists, such that 
\begin{equation}
\label{Cham}
\esp[X^s]\; =\; \frac{({\bf a})_s ({\bf b})_{-s}}{({\bf c})_s ({\bf d})_{-s}}
\end{equation}
for every $s\in (-\min({\bf a}), \min({\bf b})),$ where ${\bf a}, {\bf b}, {\bf c}$ and ${\bf d}$ are four finite and possibly empty sets of positive numbers. A challenging open problem is to characterize the existence of such distributions in terms of the four finite sets ${\bf a}, {\bf b}, {\bf c}$ and ${\bf d}.$ 

When they exist, these distributions are unique by uniqueness of the Mellin transform on the non-trivial strip $(-\min({\bf a}), \min({\bf b})).$ From the point of view of special functions, the Mellin transform is also formerly invertible in a Meijer $G-$function. More precisely, one has by definition
\begin{equation}
\label{MellG}
\int_0^\infty f(x)\, x^{s-1}\, dx\; =\;\frac{({\bf a})_s ({\bf b})_{-s}}{({\bf c})_s ({\bf d})_{-s}}
\end{equation}
for every $s\in (-\min({\bf a}), \min({\bf b})),$ where $f$ is the Meijer $G-$function
$$f(x) \; =\; G^{m, n}_{p+n, q+m} \lpa x^{-1} \lva \begin{array}{c} {\bf A} \\ {\bf B}\end{array}\right. \rpa$$
with ${\bf A} = \{1-a_n,\ldots, 1-a_1, d_1,\ldots, d_p\}$ and ${\bf B} = \{b_1,\ldots, b_m, 1-c_q,\ldots, 1-c_1\}$ finite and possibly empty sets. We refer to Chapter 16 in \cite{N} for more detail on the Meijer $G-$functions. One should remember that the above $f$ is a function strictly speaking on $(0,1)\cup(1,\infty)$ only, and that there might be a singularity at one. For example, if ${\bf a} = (1,3), {\bf c} = (2,2)$ and ${\bf b} ={\bf d} =\emptyset,$ an easy computation shows that the corresponding random variable $X$ exists with distribution
$$\pb[X\in dx]\; =\; \frac{1}{2}\, (\delta_1(dx)\, +\, \Un_{(0,1)}\, dx).$$ 
The solution to the above existence problem amounts to characterizing the non-negativity of $f$ on $(0,1)\cup(1,\infty)$ and of the singularity at one. To this end, simulations of Meijer $G-$functions can be performed by the Wolfram package \cite{Wolf}, which is however not very robust when the amount of parameters becomes large or in the presence of the singularity. The general problem of the non-negativity of Meijer functions on the positive half-line is a hard one, which can be traced back to \cite{K} for the classical hypergeometric function.

\subsection{The case $\bb =\emptyset$} Before stating, in the next paragraph, our contribution to this problem, let us give an account on the case ${\bf b} =\emptyset,$ which has been studied in \cite{CL, D, J, KP} from various points of view. It is clear that the condition $\dd =\emptyset$ is then necessary for the existence of $X.$ Moreover, considering the random variable $X^{-1},$ the problem is equivalent to the case  ${\bf a} ={\bf c} =\emptyset.$ Setting $n =\sharp\{{\bf a}\}$ and $p =\sharp\{{\bf c}\},$ we discard the obvious case $p = 0$ where $X$ always exists and is distributed as the independent product
$$X\; \elaw\; \G_{a_1}\,\times\,\cdots\, \times\, \G_{a_n}.$$
Here and throughout, $\G_t$ denotes a standard Gamma random variable with parameter $t > 0$ and we make the convention that an empty product is one. When $p\ge 1,$ the following necessary conditions for the existence of $X$ are obtained in \cite{CL}:
\begin{equation}
\label{Cle}
\lacc\begin{array}{l} p\,\le\, n,\\ a_1\, \le\, c_1, \\ a_1 +\cdots +a_p\, \le\, c_1 +\cdots +c_p.\end{array}\right.
\end{equation}
This set of conditions is easily seen to be sufficient for $p=1$ or $n=2.$ For $p = 1$ one has indeed
$$X\; \elaw\; \B_{a_1, c_1 - a_1}\,\times\, \G_{a_2}\,\times\,\cdots\, \times\, \G_{a_n}$$
where $\B_{a,b}$ denotes, here and throughout, a standard Beta random variable with parameter $a,b > 0$ and we make the conventions $\B_{0,b} = 0$ for all $b >0$ and $\B_{a,0} = 1$ for all $a > 0.$ It is interesting to mention that for $n=2,$ an example of such independent Beta-Gamma products is the square of a Brownian supremum area, which is thoroughly studied in \cite{J} - see Theorem 1.6 therein. To show the sufficiency of (\ref{Cle}) in the remaining case $p= n=2,$ we first write down the general formula
\begin{equation}
\label{Malm}
\frac{({\bf a})_s}{({\bf c})_s}\; =\; \exp -\lacc \int^{\infty}_0 (1- e^{-sx}) \lpa \frac{\varphi_\aa(x) - \varphi_\cc (x)}{x (1- e^{-x})}\rpa\, dx\racc 
\end{equation} 
for all $p=n$ with the notations
$$\varphi_\aa(x)\; =\; e^{-a_1 x}\, +\, \cdots\, +\, e^{-a_n x}\qquad\mbox{and}\qquad \varphi_\cc(x)\; =\; e^{-c_1 x}\, +\, \cdots\, +\, e^{-c_n x},$$
which follows e.g. from Theorem 1.6.2 (ii) in \cite{AAR}. By the L\'evy-Khintchine formula we see that $X$ exists, and has support $[0,1]$ as the exponential of some negative ID random variable, as soon as
\begin{equation}
\label{Schur}
\varphi_{\aa}\;\ge\;\varphi_{\cc} \quad\mbox{on $\,\rl^+.$}
\end{equation} 
This condition is then easily checked for $p = n = 2,\, a_1 \le c_1$ and $a_1 +a_2 \le c_1 + c_2.$ In the case $a_1 + a_2 = c_1 + c_2$ there is a singularity at one which is given by
$$\pb[X=1] \; =\;\lim_{s\to\infty} \esp [X^s]\; =\; \frac{\Ga(c_1)\Ga(c_2)}{\Ga(a_1)\Ga(a_2)}\; \le\; 1,$$
and the density of $X$ on $(0,1)$ can be computed in terms of a hypergeometric function $_{2}F_{1}$ - see Theorem 6.2 in \cite{D}. In the case $a_1 + a_2 < c_1 + c_2$ the random variable $X$ is absolutely continuous, with a density on $(0,1)$ which is also given in terms of a hypergeometric function $_{2}F_{1}$ - see Formula (5.1) in \cite{D}.

We next observe that the set of conditions (\ref{Cle}) is not sufficient for $p=n= 3,$ as the following counterexample shows:
$$\lacc \begin{array}{l}
a_1 = 2, a_2 = 16/5, a_3 = 17/5,\\
c_1 = 11/5, c_2 = 12/5, c_3 = 4.
\end{array}
\right.$$
Indeed, if there existed a corresponding random variable $X$, one would have
$$\pb[X=1] \; =\;\lim_{s\to\infty} \esp [X^s]\; =\; \frac{\Ga(c_1)\Ga(c_2)\Ga(c_3)}{\Ga(a_1)\Ga(a_2)\Ga(a_3)}\; = \; \frac{150}{132} \; > \; 1,$$
which implies that the density of $X$ on $(0,1)$ would have to take negative values. Even in the simpler case $p= n,$ characterizing the existence of $X$ does not seem easy in general. An interesting open question is whether the sufficient condition (\ref{Schur}) is also necessary, which would mean 
\begin{equation}
\label{ID}
X \;{\rm exists}\;\Longleftrightarrow\; \log X\; {\rm is\; ID.}
\end{equation}
In the absolutely continuous case $a_1 +\cdots + a_p < c_1 +\cdots + c_p,$ this open question is easily seen to be equivalent to Conjecture 1 in \cite{KP}. We also refer to the whole Section 2 in \cite{KP} for a set of conditions on the sets $\aa$ and $\cc$ ensuring (\ref{Schur}). The necessity of (\ref{Schur}), with the notation $\varphi_\cc(x) = e^{-c_1 x} + \cdots + e^{-c_p x},$ can also be asked for $p < n.$ Observe that the necessity of (\ref{Schur}) would also imply the equivalence (\ref{ID}) in view of the more general formula
\begin{equation}
\label{LK}
\frac{({\bf a})_s}{({\bf c})_s}\; =\; \exp \lacc (\psi(\aa) -\psi(\cc))\, s\, +\, \int_{-\infty}^0 (e^{sx} - 1 - sx) \lpa \frac{\varphi_\aa(\vert x\vert ) - \varphi_\cc (\vert x\vert )}{\vert x\vert (1- e^{-\vert x\vert})}\rpa\, dx\racc
\end{equation} 
which follows from a combination of Theorem 1.6.1 (ii) and Theorem 1.6.2 (ii) in \cite{AAR}, where $\psi$ is the digamma function and we have used the notations
$$\psi(\aa)\; =\; \psi(a_1)\, +\, \cdots\, +\, \psi(a_n)\qquad\mbox{and}\qquad \psi(\cc)\; =\; \psi(c_1)\, +\, \cdots\, +\, \psi(c_p).$$ Notice finally that for $p < n,$ the random variable $X$ has unbounded support since $(\esp[X^s])^{1/s} \to\infty$ as $s\to\infty,$ whereas its support is bounded in $[0,1]$ for $p=n$ since then $(\esp[X^s])^{1/s} \to 1.$ 
    
\subsection{A family of solutions for $\aa\neq\emptyset$ and $\bb\neq\emptyset$}
\label{A family  of solutions}
 We now come back to the general setting of Paragraph \ref{general}, and we use the notations therein. By the same arguments as in \cite{CL}, the following conditions are easily seen to be necessary:
$$
\lacc\begin{array}{l} p +q\,\le\, n+m,\\ \min(\aa)\, \le\, \min(\cc), \\ \min(\bb)\, \le\, \min (\dd).\end{array}\right.$$
Moreover, in the case $p=q=1$ these conditions are also sufficient, with
$$X\; \elaw\; \frac{\B_{a_1, c_1 - a_1}\,\times\, \G_{a_2}\,\times\,\cdots\, \times\, \G_{a_n}}{\B_{b_1, d_1 - b_1}\,\times\, \G_{b_2}\,\times\,\cdots\, \times\, \G_{b_m}}\cdot$$

\noindent
The following theorem shows however that the condition $p\le n$ which was necessary in the setting of the previous paragraph, is in general not necessary anymore.

\begin{theorem}
\label{MGT}
For every $a,b,c,d > 0,$ the distribution 
$${\rm D}\lcr\begin{array}{cc} a & b\\ (c,d) & -\end{array}\rcr\; =\; {\rm D}\lcr\begin{array}{cc} a & b\\ (d,c) & -\end{array}\rcr$$ 

\medskip

\noindent
{\em (a)} exists if $c+d\ge 3a + b + 1/2\;$ and $\;\min(c,d)\ge\min(2a +b,a+1/2).$

\medskip

\noindent
{\em (b)} does not exist if $c+d < 3a + b + 1/2\;$ or $\;\min(c,d) \le a.$

\end{theorem}

The argument for this result relies mainly on the Weber-Schafheitlin formula for the Bessel function of the first kind
$$J_\a(z)\; =\; \frac{(z/2)^\a}{\Ga(\a +1)}\,\pFq{0}{1}{-\!\! -}{\a +1}{-\lpa \frac{z^2}{2}\rpa}\; =\; \sum_{n\ge 0}  \frac{(-1)^n\, (z/2)^{2n+\a}}{n!\,\Ga(n+\a+1)}$$
defined for all $\a\in\rl$ and $z\in\CC\slash\rl^-.$ A short new proof of the latter formula is given in the Appendix for the sake of completeness. There is a direct connection between Theorem \ref{MGT} and the recent paper \cite{CY1} dealing with the non-negativity of the generalized hypergeometric function 
$$\pFq{1}{2}{a}{b,c}{x}$$ 
on the negative half-line. More precisely, we will see in (\ref{Eq1F2}) below that 
$${\rm D}\lcr\begin{array}{cc} a & b\\ (c,d) & -\end{array}\rcr\; \mbox{exists}\;\Longleftrightarrow\; \pFq{1}{2}{a+b}{c+b,d+b}{-x}\; \ge\; 0\quad\mbox{for all $x\ge 0,$}$$
and it is then easy to see that Theorem 6.1 in \cite{CY1} is equivalent to Theorem \ref{MGT}. Overall, our argument is considerably shorter than all the arguments involved in \cite{CY1}, and the above result can be viewed as a simple proof of Theorem 6.1 therein. The connection with $_{1}F_{2}$ hypergeometric functions also enables us to state the following improvement, which for its part (b) relies on the more recent and more involved paper \cite{CY2}.

\begin{theorem}
\label{MGT2}
For every $a,b > 0,$ there exists a non-increasing continuous function $f_{a,b}$ defined on $[(3a +b)/2 +1/4,\infty)$ such that 
\begin{equation}
\label{EqD}
{\rm D}\lcr\begin{array}{cc} a & b\\ (c,d) & -\end{array}\rcr\; \mbox{exists}\;\Longleftrightarrow\; f_{a,b}(d)\, \le\, c\, \le\, d \quad\mbox{or}\quad f_{a,b}(c)\,\le\, d\, \le\, c.
\end{equation}
Moreover,

\medskip

{\em (a)} for $u\in[(3a +b)/2 +1/4,\max(2a +b,a+1/2)],$ one has $f_{a,b}(u) = 3a + b + 1/2 -u.$ 

\medskip

{\em (b)} for $u >\max(2a +b,a+1/2)],$ one has $f_{a,b}(u)\; \in\; ]a, a+ (a + b)/2(u-a)].$

\medskip

\noindent
In particular, one has $f_{a,b}(u)\to a$ as $u\to\infty.$

\end{theorem}

\section{Proofs}
\subsection{Proof of Theorem 1}
\proof We begin with the existence result and we first consider the case $c = 2a + b$ and $d = a+ 1/2.$ Introducing the function
$$f(x)\; =\; \frac{\sqrt{\pi} \Ga(2a+b)\Ga(a+1/2)}{\Ga(a)\Ga(b)}\, x^{a-3/2} J_{a+b-1/2}^2(x^{-1/2})$$
which is non-negative on $(0,\infty)$ (and vanishes an infinite number of times), we get from the Weber-Schafheitlin formula -- see the Proposition in the Appendix --  and a change of variable
\begin{eqnarray*}
\int_0^\infty x^s\, f(x)\, dx & = & \frac{2 \sqrt{\pi} \Ga(2a+b)\Ga(a+1/2)}{\Ga(a)\Ga(b)}\,\int_0^\infty x^{-2a-2s} J_{a+b-1/2}^2(x)\, dx\\
& = & \frac{\Ga(2a+b)\Ga(a+1/2)\Ga(a +s)\Ga(b-s)}{\Ga(a)\Ga(b)\Ga(a+1/2 +s)\Ga(2a+b +s)}\; = \; \frac{(a)_s(b)_{-s}}{(c)_s (d)_s}
\end{eqnarray*} 
for every $s\in(-a,b).$ We next consider the case $c+d = 3a + b + 1/2$ and $\min(c,d) > \min(2a +b,a+1/2).$ It follows from the case $p=n=2,\, a_1 < c_1$ and $a_1 + a_2 = c_1 +c_2$ in the previous paragraph that the distribution
$${\rm D}\lcr\begin{array}{cc} (2a +b, a+1/2) & -\\ (c,d) & -\end{array}\rcr$$
exists and, by the previous case, so does
$${\rm D}\lcr\begin{array}{cc} a & b\\ (c,d) & -\end{array}\rcr\; =\; {\rm D}\lcr\begin{array}{cc} a & b\\ (2a +b, a+1/2) & - \end{array}\rcr\;\odot \; {\rm D}\lcr\begin{array}{cc} (2a +b, a+1/2) & -\\ (c,d) & -\end{array}\rcr,$$
where we have used the notation $\odot$, the concatenation rule and the simplification rule in \cite{CL} p.1046. The remaining case in the proof of (a) is $c+d > 3a + b + 1/2$ and $\min(c,d) \ge \min(2a +b,a+1/2),$ which is an easy consequence of the previous case since, setting $m = \min(c,d)$ and $M = \max(c,d),$ the distribution
$${\rm D}\lcr\begin{array}{cc} 3a +b +1/2 -m & -\\ M & -\end{array}\rcr$$
exists from the case $p=n=1$ and $a_1 < c_1$ in the previous paragraph, and one has  
$${\rm D}\lcr\begin{array}{cc} a & b\\ (c,d) & -\end{array}\rcr\; =\; {\rm D}\lcr\begin{array}{cc} a & b\\ (m,3a +b +1/2-m) & - \end{array}\rcr\;\odot \; {\rm D}\lcr\begin{array}{cc} 3a +b +1/2 -m & -\\ M & -\end{array}\rcr.$$

\bigskip

We now proceed to the non-existence result, which is easy for $ m =\min(c,d) \le a :$ if $m = a,$ one has
$${\rm D}\lcr\begin{array}{cc} a & b\\ (c,d) & -\end{array}\rcr\; =\; {\rm D}\lcr\begin{array}{cc} - & b\\ M & - \end{array}\rcr$$
and it is clear that such a distribution cannot exist, whereas if $m < a$ the Mellin transform
$$s\; \mapsto\; \frac{(a)_s(b)_{-s}}{(c)_s (d)_s}$$
vanishes at $s= -m$ inside the domain of convergence $(-a,b),$ which is impossible for a probability distribution. To handle the remaining case $c+d < 3a + b + 1/2,$ we shall use the following identification
$$G^{1, 1}_{1, 3} \lpa x \lva \begin{array}{c} 1-a \\ 0,1-b,1-c \end{array}\right. \rpa\; =\; \frac{\Ga(a)}{\Ga(b)\Ga(c)}\,\pFq{1}{2}{a}{b,c}{-x}$$
for every $a,b,c, x >0.$ This particular case of Formula 16.18.1 in \cite{N} will also be useful in the sequel. By (\ref{MellG}) and Formula 16.19.2 in \cite{N} with $\mu = -b,$ we deduce
\begin{equation}
\label{Eq1F2}
{\rm D}\lcr\begin{array}{cc} a & b\\ (c,d) & -\end{array}\rcr\; \mbox{exists}\;\Longleftrightarrow\; \pFq{1}{2}{a+b}{c+b,d+b}{-x}\; \ge\; 0\quad\mbox{for all $x\ge 0.$}
\end{equation}
We can now appeal to the asymptotic behaviour for the generalized hypergeometric functions, to be found e.g. in \cite{N}. Applying Formula 16.11.8 therein with $q =\kappa = 2$ and $\nu = a-b-c-d +1/2 > -2(a+b),$ we get
$$\pFq{1}{2}{a+b}{c+b,d+b}{-x}\; = \; \frac{x^{\nu/2}}{\sqrt{\pi}}\,\cos (\sqrt{x} +\nu\pi/2)\; +\; o(x^{\nu/2})\qquad\mbox{as $x\to\infty,$}$$
which shows that the function on the left-hand side takes negative values on $(0,\infty).$
\endproof

\begin{remark}{\em (a) For fixed $a,b >0,$ the existence set defined in (a) is convex with extremal points $(2a+b, a+1/2)$ and $(a+1/2, 2a +b).$ The main argument for (a) is the analysis on these extremal points, which is a consequence of the Weber-Schafheitlin formula.

\medskip

(b) As expected, the leading term in the asymptotic expansion of the hypergeometric function $_{1}F_{2}$ becomes positive when $c+d > 3a + b + 1/2$ viz. $\nu < -2(a+b).$ It can be shown from Formula 16.11.8 in \cite{N} and a finer analysis that the leading term takes negative values when $c+d = 3a + b + 1/2$ and $\min(c,d) < \min(2a +b,a+1/2).$ We omit details.

\medskip

(c) Setting $X_{a,b}$ for the random variable corresponding to the extremal distribution
$$ {\rm D}\lcr\begin{array}{cc} a & b\\ (2a +b, a+1/2) & - \end{array}\rcr,$$
the same argument leading to (\ref{LK}) implies
\begin{eqnarray*}
\esp[X_{a,b}^s] & = & \exp \lacc (\psi(a) -\psi(b,c,d))\, s\, +\, \int_{-\infty}^0 (e^{sx} - 1 - sx) \lpa \frac{\varphi_a(\vert x\vert) - \varphi_{c,d} (\vert x\vert)}{\vert x\vert (1- e^{-\vert x\vert})}\rpa\, dx\right.\\
& & \qquad\qquad\qquad \, +\, \left. \int^{\infty}_0 (e^{sx} - 1 - sx) \lpa \frac{\varphi_b(x)}{x (1- e^{-x})}\rpa\, dx\racc
\end{eqnarray*}
for every $s\in (-a,b).$ In the recent terminology of \cite{LPS}, this means that $\log X_{a,b}$ is quasi-infinitely divisible (QID) with quasi-L\'evy measure having density
$$\frac{1}{\vert x\vert (1- e^{-\vert x\vert})} \lpa (\varphi_a(\vert x\vert) - \varphi_{c,d} (\vert x\vert)) \Un_{\{x<0\}}\, +\,  \varphi_b(x)\Un_{\{x >0\}}\rpa.$$
This density, which is negative in an interval $(a_*, 0)$ with $a_* <0,$ vanishes at $a_*,$ and is positive otherwise, tends to $-\infty$ and is not integrable at $0\!-\!.$ This implies that $\log X_{a,b}$ is not ID and that (\ref{ID}) fails. It would be interesting to construct more general families of QID distributions with the help of the Gamma function.}

\end{remark}
\subsection{Proof of Theorem 2}

\proof
Fix $a,b > 0$ and consider
$$\cD_{a,b}\; =\; \lacc(c,d)\in (a,\infty)\times (a,\infty)\; \mbox{such that}\;{\rm D}\lcr\begin{array}{cc} a & b\\ (c,d) & -\end{array}\rcr\; \mbox{exists}\racc.$$
It is clear from the definition that the set $\cD_{a,b}$ is closed and symmetric with respect to the line $\{c=d\}.$ Theorem 1 shows that $\cD_{t,a,b} =\cD_{a,b}\cap\{c+d =t\}$ is non-empty if and only if $t\ge 3a+b + 1/2.$ The argument for the second case in Part (a) of Theorem 1 shows clearly that $\cD_{t,a,b}$ is a closed segment $[(x_t,y_t), (y_t,x_t)]$ with $x_t \in (a,m]$ for every $t\ge 3a +b +1/2.$ Besides, the easy fact that 
$$(c,d)\,\in\,\cD_{a,b}\;\Longrightarrow\; (c,d+s)\,\in\,\cD_{a,b}\;\;\mbox{for every $s\ge 0$}$$
which was used in the argument for the third case in Part (a) of Theorem 1,
implies that the function $t\mapsto x_t$ is non-increasing and that the function $t\mapsto y_t$ is increasing. Finally, both functions are clearly continuous by the closedness of $\cD_{a,b}.$ Consider now the mapping
$$f_{a,b}(u)\; =\; \lacc\begin{array}{ll} 3a + b + 1/2 -u & \mbox{if $u\in[(3a +b)/2 +1/4,\max(2a +b,a+1/2)]$}\\ 
x(y^{-1}(u))  & \mbox{if $u >\max(2a +b,a+1/2),$}
\end{array}\right.$$
which is continuous and non-increasing on $[(3a +b)/2 +1/4,\infty).$ Putting everything together with Theorem 1 implies the required equivalence (\ref{EqD}), and it is clear that (a) is fulfilled. To show (b) and conclude the proof of the Theorem, it suffices to combine (\ref{Eq1F2}) and Theorem 4.2 in \cite{CY2}.
 
\endproof

From the above discussion, the latter result can also be expressed in terms of the hypergeometric function $_{1}F_{2}.$ More precisely, it follows from (\ref{Eq1F2}) that 
$$\pFq{1}{2}{2a}{b,c}{-x}\; \ge\; 0\quad\mbox{for all $x\ge 0$}\;\Longleftrightarrow\; (b-a,c-a)\;\in \;\cD_{a,a}$$
for every $a,b,c > 0.$ This functional representation for the diagram of non-negativity of $_{1}F_{2}$ on the negative half-line seems unnoticed. Following the introduction in \cite{CY2}, one can also rephrase a problem of Askey and Szeg\"o  in the following way : for every $a, b >0$ one has
$$\int_0^x t^{b-a} J_{a+b-1}(t)\, dt\; \ge\; 0\quad\mbox{for all $x\ge 0$}\;\;\Longleftrightarrow\;\; (a,1)\;\in \;\cD_{b,b}.$$
Perhaps can this probabilistic reformulation of an old and celebrated problem on Bessel functions of the first kind - see also Chapter 7.6 in \cite{AAR} for a discussion - be useful. We would like to finish this paper with the following open problem, which is natural in view of Remark (a).

\begin{conjecture}
For every $a,b >0,$ the set $\cD_{a,b}$ is convex. 
\end{conjecture}

This conjecture amounts to the convexity of the function $f_{a,b}.$ This would imply that the latter function is also decreasing, in other words that the non-increasing function $t\mapsto x_t$ which was introduced during the proof of Theorem 2 is actually decreasing. But we were not able to prove this.
 
\section*{Appendix}

In this Appendix, we give an independent proof of the classical Weber-Schafheitlin formula, which is the exact computation of a certain Mellin transform. Let us mention that this formula has been given several generalizations over the years - see e.g. \cite{B,SE,KR} and the references therein. The latter, however, do not seem to have direct applications in our setting.
 
\begin{proposition}[Weber-Schafheitlin] 
For every $\a > -1/2$ and $s\in (0, \a +1/2),$ one has
$$\int_0^\infty z^{-2s}\, J_\a^2(z)\, dz\; =\; \frac{\Ga(s)\Ga(\a+1/2 -s)}{2\sqrt{\pi}\Ga(1/2 +s) \Ga(\a+1/2 +s)}\cdot$$ 
\end{proposition}

The proof relies on the following  computation of a Fourier transform, which is originally due to von Lommel. Observe in passing that this simple computation provides a family of solutions to the so-called Van Dantzig problem - see Theorem 1 and Corollary 2 in \cite{Ju}.

\begin{lemma}[von Lommel]
\label{Lomme}
 For every $\a > -1/2$ and $z\in \CC\slash\rl^-,$ one has
$$\frac{1}{\sqrt{\pi}\,\Ga(\a +1/2)}\, \lpa\frac{z}{2}\rpa^\a \int_{-1}^1 e^{\ii tz} (1-t^2)^{\a -1/2}\, dt\; =\; J_\a(z) .$$ 
\end{lemma}

\proof Expanding $e^{\ii tz}$ as a series, switching the sum and the integral, cancelling the odd terms and changing the variable, the right-hand side transforms into
\begin{eqnarray*}
\frac{1}{\sqrt{\pi}\,\Ga(\a +1/2)} \lpa\frac{z}{2}\rpa^\a \sum_{n\ge 0} (-1)^n \frac{z^{2n}}{(2n)!}\, \int_0^1 t^{n-1/2}(1-t)^{\a -1/2}\, dt & = & \!\lpa\frac{z}{2}\rpa^\a\, \sum_{n\ge 0} (-1)^n \frac{\Ga(n+1/2)\,z^{2n}}{\sqrt{\pi}\,(2n)!\,\Ga(n+\a+1)}\\
& = & \sum_{n\ge 0}  \frac{(-1)^n\, (z/2)^{2n+\a}}{n!\,\Ga(n+\a+1)}\; =\; J_\a(z),
\end{eqnarray*}
where the second equality follows from Gauss' multiplication formula - see e.g. Theorem 1.5.1 in \cite{AAR}.

\endproof

We now proceed to our proof of the Proposition. The original argument is the consequence of a more general result involving a quadratic transformation of some hypergeometric function - see \cite{Wot} p.402 for details and also Exercise 4.15 in \cite{AAR} for a more modern presentation of this result. We follow here another, apparently unnoticed and overall simpler argument relying on the above Lommel formula, the Fresnel integral and the Selberg integral. For every $\a > -1/2$ and $s\in (0, \a +1/2),$ the above Lemma implies first 
$$\int_0^\infty z^{-2s}\, J_\a^2(z)\, dz\; =\; \frac{1}{\pi 4^\a \Ga(\a +1/2)^2}\int_0^\infty z^{2\a - 2s}\lpa \int_{[-1,1]^2}\!\! \cos(z(t+u)) ((1-t^2)(1-u^2))^{\a -1/2}dt\, du\rpa  dz.$$    
Supposing next $s > \a,$ this transforms into
\begin{eqnarray*}
\int_0^\infty z^{-2s}\, J_\a^2(z)\, dz & = & \frac{1}{\pi 4^\a \Ga(\a +1/2)^2} \int_{[-1,1]^2}\!\! ((1-t^2)(1-u^2))^{\a -1/2} \lpa \int_0^\infty z^{2\a - 2s} \cos(z(t+u)) \, dz\rpa dt\, du\\
& = & \frac{1}{4^{\a+1/2} \Ga(\a +1/2)^2 \Ga(2s-2\a)\cos(\pi(s-\a))}\\
& & \qquad\qquad\qquad\qquad\qquad\qquad \times\; \int_{[-1,1]^2} ((1-t^2)(1-u^2))^{\a -1/2} \vert t-u\vert^{2s -2\a-1} dt\, du
\end{eqnarray*}
where in the second equality we have used the Fresnel integral - see e.g. Exercise 1.19 in \cite{AAR} - and the change of variable $u\mapsto -u.$ On the other hand, one has
$$4^{\a+1/2} \Ga(2s-2\a)\cos(\pi(s-\a))\; =\; \frac{4^{\a+1/2} \Ga(2s-2\a)}{\Ga(s-\a+1/2)\Ga(\a-s +1/2)}\; =\; \frac{4^s\sqrt{\pi}\, \Ga(s-\a)}{\Ga(\a +1/2-s)}$$
by the complement and multiplication formulas, for the Gamma function. Therefore,
\begin{eqnarray*}
\int_0^\infty z^{-2s}\, J_\a^2(z)\, dz & = &  \frac{\Ga(\a +1/2-s)}{4^s\sqrt{\pi}\,\Ga(\a +1/2)^2  \Ga(s-\a)} \int_{[-1,1]^2} ((1-t^2)(1-u^2))^{\a -1/2} \vert t-u\vert^{2s -2\a-1} dt\, du\\
 & = &  \frac{4^{\a-1/2}\Ga(\a +1/2-s)}{\sqrt{\pi}\,\Ga(\a +1/2)^2  \Ga(s-\a)} \int_{[0,1]^2} (t(1-t)u(1-u))^{\a -1/2} \vert t-u\vert^{2s -2\a-1} dt\, du\\
& = &  \frac{4^{\a}\Ga(s)^2\Ga(2s-2\a)\Ga(\a+1/2 -s)}{2\sqrt{\pi}\,\Ga(2s)\Ga(s-\a)\Ga(s-\a+1/2)\Ga(\a+1/2 +s)}\\
& = &  \frac{\Ga(s)\Ga(\a+1/2 -s)}{2\sqrt{\pi}\,\Ga(1/2 +s) \Ga(\a+1/2 +s)}
\end{eqnarray*} 
for every $s\in(\max(0,\a), \a+1/2),$ where in the third equality we have used the Selberg integral - see e.g. Theorem 8.1.1 in \cite{AAR} - and the fourth equality follows again from Gauss' multiplication formula. By analytic continuation, the formula remains true for every $s\in(0, \a+1/2).$

\qed

\end{document}